\documentstyle[12pt]{article}
\makeatletter
\renewcommand\thesection{\@Roman\c@section}
\renewcommand\thesubsection{\thesection.\@arabic\c@subsection}
\makeatother

\begin{document}
\begin{titlepage}
\begin{flushright}
math.QA/9908008
\end{flushright}
\vskip.3in

\begin{center}
{\Large \bf Vertex Operators of $U_q[\widehat{gl(N|N)}]$ and Highest
Weight Representations of $U_q[\widehat{gl(2|2)}]$ }
\vskip.3in
{\large Wen-Li Yang $^{a,b,c}$ and Yao-Zhong Zhang $^{c}$}
\vskip.2in
{\em $~^a$ CCAST (World Laboratory), P.O.Box 8730, Beijing 100080, China\\
$~^b$ Institute of Modern Physics, Northwest University, 
Xian 710069 ,China \\
$~^c$ Department of Mathematics, University of Queensland, Brisbane,
     Qld 4072, Australia}
\end{center}

\vskip 2cm
\begin{center}
{\bf Abstract}
\end{center}
We determine the exchange relations of the level-one q-vertex operators of
the quantum affine superalgebra $U_q[\widehat{gl(N|N)}]$. We study in
details the level-one irreducible highest weight representations of 
$U_q[\widehat{gl(2|2)}]$, and compute the characters and supercharacters
associated with these irreducible modules.


\end{titlepage}


\def\a{\alpha}
\def\b{\beta}
\def\d{\delta}
\def\e{\epsilon}
\def\ve{\varepsilon}
\def\g{\gamma}
\def\k{\kappa}
\def\l{\lambda}
\def\o{\omega}
\def\t{\theta}
\def\s{\sigma}
\def\D{\Delta}
\def\L{\Lambda}

\def\R{\overline{R}}
\def\G{{gl(N|N)}}
\def\hS{{\widehat{sl(N|N)}}}
\def\hG{{\widehat{gl(N|N)}}}
\def\R{{\cal R}}
\def\hR{{\hat{\cal R}}}
\def\C{{\bf C}}
\def\P{{\bf P}}
\def\Z2{{{\bf Z}_2}}
\def\Z{{\bf Z}}
\def\T{{\cal T}}
\def\H{{\cal H}}
\def\F{{\cal F}}
\def\V{\overline{V}}
\def\trho{{\tilde{\rho}}}
\def\tphi{{\tilde{\phi}}}
\def\tT{{\tilde{\cal T}}}
\def\uqsnh{{U_q[\widehat{sl(N|N)}]}}
\def\uqgnh{{U_q[\widehat{gl(N|N)}]}}
\def\uq1h{{U_q[\widehat{gl(1|1)}]}}
\def\uqg2h{{U_q[\widehat{gl(2|2)}]}}


\def\beq{\begin{equation}}
\def\eeq{\end{equation}}
\def\bea{\begin{eqnarray}}
\def\eea{\end{eqnarray}}
\def\ba{\begin{array}}
\def\ea{\end{array}}
\def\no{\nonumber}
\def\lt{\left}
\def\rt{\right}
\newcommand{\bq}{\begin{quote}}
\newcommand{\eq}{\end{quote}}

\newtheorem{Theorem}{Theorem}
\newtheorem{Definition}{Definition}
\newtheorem{Proposition}{Proposition}
\newtheorem{Lemma}{Lemma}
\newtheorem{Corollary}{Corollary}
\newcommand{\proof}[1]{{\bf Proof. }
        #1\begin{flushright}$\Box$\end{flushright}}

\newcommand{\sect}[1]{\setcounter{equation}{0}\section{#1}}
\renewcommand{\theequation}{\thesection.\arabic{equation}}

\sect{Introduction\label{intro}}
\subsection{Free bosonic representations}
Infinite-dimensional highest weight representations and the 
corresponding  vertex operators \cite{Fre92} of quantum affine
(super)algberas are two ingredients of great importance in the 
algebraic analysis of lattice integrable models \cite{Dav93,Jim94}. 
Under some assumptions on the physical space of states of the models, 
this algebraic analysis method enables one to compute the correlation 
functions and form factors of these models in the form of 
integrals \cite{Dav93,Jim94,Koy94,Yan99,Hou99}.

A powerful approach for studying the highest weight representations 
and vertex operators is the bosonization technique
\cite{Fre88,Ber89} which allows one to explicitly construct these 
objects in terms of the q-deformed free bosonic fields.

Free bosonic realization of level-one representations and the associated 
vertex operators  has been constructed for most (twisted or untwisted) 
quantum affine algebras \cite{Fre88,Ber89,Awa94,Jin95,Jin97}. This 
bosonization technique has recently been extended to the type I quantum 
affine superalgebras $U_q[\widehat{sl(M|N)}]$, $M\neq N$ \cite{Kim97} 
and $\uqgnh$ \cite{Zha98}. Level-one representations and vertex 
operators of these two cases have been constructed in terms of 
q-deformed free bosonic fields. (See also \cite{Awa97} and 
\cite{Gad97} for the  special $U_q[\widehat{sl(2|1)}]$ 
and $U_q[\widehat{gl(2|2)}]$ cases , respectively.) 

\subsection{Highest weight irreducible representations}
For most quantum affine bosonic algebras, the level-one 
free bosonic representations are already irreducible. However, for 
$U_q[\widehat{sl(M|N)}]$, $M\neq N$ or $M=N$, such 
free bosonic representations, which were constructed in \cite{Kim97,Zha98}, are 
not irreducible in general. It is interesting (and important for 
applications) to construct irreducible highest weight 
representations out of the reducible ones. 
In \cite{Kim97,Yan99} and 
\cite{Yan99b} the level-one irreducible highest weight 
representations of $U_q[\widehat{sl(2|1)}]$ and 
$U_q[\widehat{gl(1|1)}]$  have been investigated in details. 
The method adopted in these papers is to look for suitable  projection 
operators and then applying them to decompose the reducible
representations into a direct sum of irreducible ones. The characters 
and supercharacters associated with these irreducible representations 
can then be computed by the BRST resolution.

\subsection{The present work}
In the first part of this work, we determine the exchange relations 
of level-one q-vertex operators of $\uqgnh$ by using the free bosonic 
realization given in \cite{Zha98}. We also give a Miki's construction 
of the super Reshetikhin-Semenov-Tian-Shansky (RS) algebra 
\cite{Zha97} for $\uqgnh$ at level-one. The second part of this work 
is devoted to a detailed study of level-one irreducible 
highest weight representations of $\uqg2h$ and the corresponding 
vertex operators. We compute the characters and supercharacters 
associated with all these irreducible representations. 

The paper is organized as follow. In section II we briefly review
the bosonization of $\uqgnh$ given in \cite{Zha98}. In section III
we obtain the free bosonic realization of the level-one vertex operators
of $\uqgnh$ by extending the results in \cite{Zha98} for
$U_q[\widehat{sl(N|N)}]$. In section IV we determine the exchange
relations  of these vertex operators. In section V, we study 
in details the level-one irreducible highest weight representations of 
$\uqg2h$  and the associated with vertex operators. The
characters and supercharacters corresponding to these irreducible
modules are computed.

\sect{ Bosonization of $\uqgnh$ at Level-One}
\subsection{Drinfeld basis of  $\uqgnh$}
Let $\{\a_i,~i=0,1,\cdots,2N-1\}$ denote a
chosen set of simple roots of the affine superalgebra
$\hS$ \cite{Kac78}.  As in \cite{Zha98}, 
we choose the following non-standard  system of simple
roots, all of which are odd, 
\bea
&&\a_0=\d-\ve_1+\ve_{2N}\,,\no\\
&&\a_l=\ve_l-\ve_{l+1}\,,~~~~~l=1, 2,\cdots,2N-1\label{roots}
\eea
with $\d,~\{\ve_k\}_{k=1}^{2N}$ satisfying $
(\d,\d)=(\d,\ve_k)=0$, $(\ve_k,\ve_{k'})=(-1)^{k+1}\d_{kk'}$. 
Exending the Cartan subalgebra of $\widehat{sl(N|N)}$ by adding to it 
the element 
\beq
\a_{2N}=\sum_{k=1}^{2N}\ve_k\label{root-2n}~,
\eeq
we obtain an enlarged  Cartan matrix with elements  
$a_{ij}=(\a_i,\a_j)$, $i,j=0,1,2,\cdots,2N$. Explicitly, 
\bea
(a_{ij})=\left(
\begin{array}{ccccccccc}  
0 & -1 &  &  & &1&-2&\\
-1&  0 &1 &  & && 2&\\
 &1 & 0 & \ddots & &&-2& \\
 && \ddots & \ddots & \ddots &  & & \\
 &&        &  1     &       0&-1&-2& \\
1&&        &        &      -1& 0&2&\\
-2&2&-2    &\cdots  &\cdots  & 2&0&
\end{array}
\right) ~~.\label{Cartan-matrix}
\eea
\noindent Then the Cartan matrix $(a_{ij})$, $i,j=1,2,\cdots,2N$ is
invertible, which 
is essential in the construction of q-vertex operators.

We denote by $\H$ the extended Cartan
subalgebra and by $\H^*$ the dual of $\H$.
Let $\{h_0,h_1,\cdots,h_{2N},d\}$ be a basis of $\H$, where 
$h_{2N}$ is the element in $\H$ corresponding to $\a_{2N}$ and $d$
is the usual derivation operator.
We identify the dual $\H^*$ with $\H$ via the form $(~,~)$.  
Let $\{\L_0,\L_1,\cdots,\L_{2N},\delta\}$
be the dual basis with $\L_j$ being fundamental weights. Explicitly
\cite{Zha98}
\bea
&&\L_{2N}=\frac{1}{2N}\sum_{k=1}^{2N}(-1)^{k+1}\ve_k,\no\\
&&\L_i=\L_0+\sum_{k=1}^i(-1)^{k+1}\ve_k-\frac{i}{2N}\sum_{k=1}^{2N}
   (-1)^{k+1}\ve_k,
\eea
where $i=0,1,\cdots,2N-1$. Obviously,
$\H$ ($\H^*$) constitutes  the (dual) Cartan
subalgebra of $\hG$.

The quantum affine superalgebra $\uqgnh$ is a quantum (or $q$-)
deformation of the universal enveloping algebra of $\hG$ and is
generated by the Chevalley generators $\{e_i,\;f_i\;q^{h_j},\;d|
i=0,1,\cdots,2N-1,\;j=0,1,\cdots,2N\}$.
The $\Z_2$-grading of the Chevalley generators is $[e_i]=[f_i]=1,~
i=0,1,\cdots,2N-1$ and zero
otherwise. The defining relations are
\bea
&&hh'=h'h,~~~~~~\forall h\in \H,\no\\
&&q^{h_j}e_iq^{-h_j}=q^{a_{ij}}e_i,~~~~[d, e_i]=\d_{i0}e_i,\no\\
&&q^{h_j}f_iq^{-h_j}=q^{-a_{ij}}f_i,~~~~[d,f_i]=-\d_{i0}f_i,\no\\
&&[e_i,f_{i'}]=\d_{ii'}\frac{q^{h_i}-q^{-h_i}}{q-q^{-1}},\no\\
&&[e_i,e_{i'}]=[f_i,f_{i'}]=0,~~~~{\rm for}~~a_{ii'}=0,\no\\
&&[[e_0,e_1]_{q^{-1}}, [e_0,e_{2N-1}]_q]=0,\no\\
&&[[e_l,e_{l-1}]_{q^{(-1)^l}}, [e_l,e_{l+1}]_{q^{(-1)^{l+1}}}]=0,\no\\
&&[[e_{2N-1},e_{2N-2}]_{q^{-1}}, [e_{2N-1},e_0]_q]=0,\no\\
&&[[f_0,f_1]_{q^{-1}}, [f_0,f_{2N-1}]_q]=0,\no\\
&&[[f_l,f_{l-1}]_{q^{(-1)^l}}, [f_l,f_{l+1}]_{q^{(-1)^{l+1}}}]=0,\no\\
&&[[f_{2N-1},f_{2N-2}]_{q^{-1}}, [f_{2N-1},f_0]_q]=0,~~~~l=1,2,\cdots,2N-2.
\eea
Here and throughout, $[a,b]_x\equiv ab-(-1)^{[a][b]}x ba$ 
and $[a,b]\equiv [a,b]_1$. 

$\uqgnh$ is a $\Z_2$-graded quasi-triangular Hopf algebra endowed with
the following coproduct $\D$, counit $\e$ and antipode $S$:
\bea
\D(h)&=&h\otimes 1+1\otimes h,\no\\
\D(e_i)&=&e_i\otimes 1+q^{h_i}\otimes e_i,~~~~
\D(f_i)=f_i\otimes q^{-h_i}+1\otimes f_i,\no\\
\e(e_i)&=&\e(f_i)=\e(h)=0,\no\\
S(e_i)&=&-q^{-h_i}e_i,~~~~S(f_i)=-f_iq^{h_i},~~~~S(h)=-h,\label{e-s} 
\eea
where $i=0,1,\cdots,2N-1$ and $h\in \H$.
Notice that the antipode $S$ is a $\Z_2$-graded algebra anti-homomorphism.
Namely, for any homogeneous elements $a,b\in\uqgnh$
$S(ab)=(-1)^{[a][b]}S(b)S(a)$, which extends to inhomogeneous elements
through linearity.
The multiplication rule on the tensor products is $\Z_2$-graded:
$(a\otimes b)(a'\otimes b')=(-1)^{[b][a']}(aa'\otimes bb')$ for any
homogeneous elements $a,b,a',b'\in \uqgnh$

$\uqgnh$ can also be realized in terms of the Drinfeld generators
\cite{Dri88} $\{X^{\pm,i}_m,\; H^j_n,\;$ $q^{\pm H^j_0}$, $c,\; d |
m\in {\bf Z},\; n\in{\bf Z}-\{0\},\; i=1,2,\cdots,2N-1,\;
j=1,2,\cdots,2N\}$.  The $\Z_2$-grading of the Drinfeld generators is
given by $[X^{\pm,i}_m]=1$ for all $i=1,\cdots,2N-1,\;m\in{\bf Z}$ and 
$[H^j_n]=[H^j_0]=[c]=[d]=0$ for all $j=1,\cdots,2N,\; n\in{\bf Z}-\{0\}$.
The relations 
satisfied by the Drinfeld generators read \cite{Yam96,Zha97,Zha98}
\bea
&&[c,a]=[d,H^j_0]=[H^j_0, H^{j'}_n]=0,~~~~\forall a\in\uqgnh\no\\
&&q^{H^j_0}X^{\pm, i}_nq^{-H^j_0}=q^{\pm a_{ij}}X^{\pm, i}_n,\no\\
&&[d,X^{\pm,i}_n]=nX^{\pm,i}_n,~~~[d,H^j_n]=nH^j_n,\no\\
&&[H^j_n, H^{j'}_m]=\d_{n+m,0}\frac{[a_{jj'}n]_q[nc]_q}{n},\no\\
&&[H^j_n,
   X^{\pm, i}_m]=\pm\frac{[a_{ij}n]_q}{n}X^{\pm,i}_{n+m}q^{\mp|n|c/2},\no\\
&&[X^{+,i}_n, X^{-,i'}_m]=\frac{\d_{ii'}}{q-q^{-1}}\lt(q^{\frac{c}{2}(n-m)}
  \psi^{+,i}_{n+m}-q^{-\frac{c}{2}(n-m)}\psi^{-,i}_{n+m}\rt),\no\\
&&[X^{\pm,i}_n, X^{\pm, i'}_m]=0,~~~~{\rm for}~a_{ii'}=0,\no\\
&&[X^{\pm,i}_{n+1}, X^{\pm,i'}_m]_{q^{\pm a_{ii'}}}
  -[X^{\pm,i'}_{m+1}, X^{\pm,i}_n]_{q^{\pm a_{ii'}}}=0,\no\\
&&[[X^{\pm,l}_m,X^{\pm,l-1}_{m'}]_{q^{(-1)^l}},[X^{\pm,l}_n,
  X^{\pm,l+1}_{n'}]_{q^{(-1)^{l+1}}}]\no\\
&&~~~+[[X^{\pm,l}_n,X^{\pm,l-1}_{m'}]_{q^{(-1)^l}},[X^{\pm,l}_m,
  X^{\pm,l+1}_{n'}]_{q^{(-1)^{l+1}}}]=0, ~~~l=2,\cdots,2N-2.
  \label{drinfeld}
\eea
where $[x]_q=(q^x-q^{-x})/(q-q^{-1})$ and
 $\psi^{\pm,j}_{n}$ are related to $H^j_{\pm n}$ by relations
\beq
\sum_{n\in{\bf Z}}\psi^{\pm,j}_{n}z^{- n}=q^{\pm H^j_0}\exp\lt(
  \pm(q-q^{-1})\sum_{n>0}H^j_{\pm n}z^{\mp n}\rt).
\eeq

The Chevalley generators are related to the Drinfeld generators by the
formulae
\bea
&&h_i=H_0^i,~~~e_i=X^{+,i}_0,~~~f_i=X^{-,i}_0,~~~
  i=1,2,\cdots,2N-1,\\
&&h_{2N}=H^{2N}_0,~~~h_0=c-\sum_{k=1}^{2N-1}H^k_0,\no\\
&&e_0=[X^{-,2N-1}_0,[X^{-,2N-2}_0,\cdots,[X^{-,3}_0,[X^{-,2}_0,
  X^{-,1}_1]_{q}]_{q^{-1}}\cdots]_{q}]_{q^{-1}}
  \,q^{-\sum_{k=1}^{2N-1}H^k_0},\no\\
&&f_0=(-1)^Nq^{\sum_{k=1}^{2N-1}H^k_0}[[\cdots [[X^{+,1}_{-1},
  X^{+,2}_0]_{q^{-1}},
  X^{+,3}_0]_{q},\cdots,X^{+,2N-2}_0]_{q^{-1}},X^{+,2N-1}_0]_{q}.\no\\
\eea
\subsection{Free Bosonic realization of  
$\uqgnh$ at level one}
In this subsection, we briefly review the bosonization of $\uqgnh$ at
level one \cite{Zha98}.

Introduce bosonic oscillators $\{a^j_n,\;c^l_n,\;Q_{a^j},\;
Q_{c^l}|n\in{\bf Z}, j=1,2,\cdots, 2N,\;l=1,2,\cdots,N\}$ which
satisfy the commutation relations
\bea
&&[a^j_n, a^{j'}_m]=(-1)^{j+1}\d_{jj'}\d_{m+n,0}\frac{[n]^2_q}{n},~~~~~
  [a^j_0, Q_{a^{j'}}]=\d_{jj'},\\
&&[c^l_n, c^{l'}_m]=\d_{ll'}\d_{n+m,0}\frac{[n]_q^2}{n},~~~~~
  [c^l_0, Q_{c^{l'}}]=\d_{ll'}.\label{oscilators}
\eea
The remaining commutation relations are zero. Associated with these  
bosonic oscillators are   the  q-deformed free bosonic fields 
\bea
&&H^j(z;\k)=Q_{A^j}+A^j_0\ln z   
   -\sum_{n\neq 0}\frac{A^j_n}{[n]_q}q^{\k |n|}z^{-n},\\
&&c^l(z)=Q_{c^l}+c^l_0\ln z-\sum_{n\neq 0}\frac{c^l_n}{[n]_q}  
   z^{-n},\\
&&H^j_\pm(z)=\pm(q-q^{-1})\sum_{n>0}A^j_{\pm n}z^{\mp n}\pm A^i_0\ln q.
\eea
\noindent Here, 
\bea
&&A^i_n=(-1)^{i+1}\lt(a^i_n+a^{i+1}_n\rt),~~~~
  Q_{A^i}=Q_{a^i}-Q_{a^{i+1}},~~~i=1,2,\cdots,2N-1,\\
&&A^{2N}_n=\frac{q^n+q^{-n}}{2}\sum_{l=1}^{2N} (-1)^{l+1}a^l_n,~~~~
  Q_{A^{2N}}=\sum_{l=1}^{2N}Q_{a^l}.\label{A-a}
\eea
\noindent Let us  define  the Drinfeld currents
\begin{eqnarray*}
&&X^{\pm,i}(z)=\sum_{n\in\Z}X^{\pm,i}_nz^{-n-1}~,~~i=1,2,\cdots,2N-1,\\
&&\psi^{\pm,j}(z)=\sum_{n\in\Z}\psi^{\pm,j}_nz^{-n}~,~~j=1,2,\cdots,2N.
\end{eqnarray*}
Then , 
\begin{Theorem}\label{free-boson}(\cite{Zha98}):
The Drinfeld generators of $\uqgnh$  at level one are
realized by the free boson fields as
\bea
&&c=1,\no\\
&&\psi^{\pm,j}(z)=e^{H^j_\pm(z)},~~~~j=1,2,\cdots,2N,\no\\
&&X^{\pm,i}(z)=:e^{\pm H^i(z;\mp\frac{1}{2})}\;Y^{\pm,i}(z):F^{\pm,i},~~~~
   i=1,2,\cdots,2N-1,
\eea
where
\bea
&&F^{\pm,2k-1}=\prod_{l=1}^{k-1}e^{\pm\sqrt{-1}\pi a^{2l-1}_0},~~~~~
   F^{\pm,2k}=\prod_{l=1}^{k}e^{\mp\sqrt{-1}\pi a^{2l-1}_0},\no\\
&&Y^{+,2k-1}(z)=e^{c^k(z)},\no\\
&&Y^{-,2k-1}(z)=\frac{1}{z(q-q^{-1})}\lt(e^{-c^k(qz)}-
  e^{-c^k(q^{-1}z)}\rt),\no\\
&&Y^{+,2k}(z)=Y^{-,2k-1}(z)=\frac{1}{z(q-q^{-1})}\lt(e^{-c^k(qz)}-
  e^{-c^k(q^{-1}z)}\rt),\no\\
&&Y^{-,2k}(z)=-Y^{+,2k-1}(z)=-e^{c^k(z)},~~~~k=1,2,\cdots,N.
   \label{boson}
\eea
\end{Theorem}

\sect{Bosonization of Level-One Vertex Operators}

We consider the evaluation representation $V_z$ of $\uqgnh$,
where $V$ is an $2N$-dimensional graded
vector space with basis vectors $\{v_1,v_2,\cdots,v_{2N}\}$. The
$\Z_2$-grading of the basis vectors is chosen to be $[v_j]=\frac{(-1)^j
+1}{2}$. Let $e_{j,j'}$ be the $2N\times 2N$ matrices satisfying
$e_{i,j}v_k=\d_{jk}v_i$.
Denote by $V^{*S}$  the left dual module of $V$ defined by
\beq
(a\cdot v^*)(v)=(-1)^{[a][v^*]}v^*(S(a)v),~~~\forall a\in \uqgnh,\;
   v\in V,\; v^*\in V^*.
\eeq
Namely, the representations on
$V^{*S}$  are given by
\beq
\pi_{V^{*S}}(a)=\pi_V(S(a))^{st},~~~~\forall a\in\uqgnh,
\eeq
where $st$ denotes the supertansposition defined by
$(A_{i,j})^{st}=(-1)^{[j]([i]+[j])}A_{j,i}$. Let $V^{*S}_z$
be the $2N$-dimensional evaluation module corresponding to $V^{*S}$.

In the homogeneous gradation, 
the Drinfeld generators are represented \cite{Zha98}
on $V_z$ by
\bea
H^i_m&=&(-1)^{i+1}\frac{[m]_q}{m}q^{(-1)^im}\lt(q^{x_i}z\rt)^m
  (e_{i,i}+e_{i+1,i+1}),\no\\
H^{2N}_m&=&z^m\frac{[2m]_q}{m}\lt[-q^m\sum_{l=1}^Ne_{2l,2l}\rt.\no\\
& &\lt. +\sum_{l=1}^N\lt(1-N+(l-1)(1-q^m)\rt)(e_{2l-1,2l-1}+
  e_{2l,2l})\rt],\no\\
H^i_0&=&(-1)^{i+1}(e_{i,i}+e_{i+1,i+1}),~~~~H^{2N}_0=\sum_{k=1}^{2N}
  (-1)^{k+1}e_{k,k},\no\\
X^{+,i}_m&=&\lt(q^{x_i}z\rt)^me_{i,i+1},~~~~
  X^{-,i}_m=(-1)^{i+1}\lt(q^{x_i}z\rt)^me_{i+1,i},\no
\eea
and on $V^{*S}_z$ by
\bea
H^i_m&=&(-1)^{i}\frac{[m]_q}{m}q^{(-1)^{i+1}m}\lt(q^{-x_i}z\rt)^m
  (e_{i,i}+e_{i+1,i+1}),\no\\
H^{2N}_m&=&-z^m\frac{[2m]_q}{m}\lt[-q^{-m}\sum_{l=1}^Ne_{2l,2l}\rt.\no\\
& &\lt. +\sum_{l=1}^N\lt(1-N+(l-1)(1-q^{-m})\rt)(e_{2l-1,2l-1}+
  e_{2l,2l})\rt],\no\\
H^i_0&=&(-1)^{i}(e_{i,i}+e_{i+1,i+1}),~~~~H^{2N}_0=\sum_{k=1}^{2N}
  (-1)^{k}e_{k,k},\no\\
X^{+,i}_m&=&-(-1)^iq^{(-1)^i}\lt(q^{-x_i}z\rt)^me_{i+1,i},~~~~
  X^{-,i}_m=-q^{(-1)^{i+1}}\lt(q^{-x_i}z\rt)^me_{i,i+1},\no
\eea
where $i=1,\cdots,2N-1$,  
$x_i=\sum_{l=1}^i (-1)^{l+1}=\frac{(-1)^{i+1}+1}{2}$.

Let
$V(\l)$ be the highest weight $\uqgnh$-module with the highest weight
$\l$. Consider the following intertwiners of $\uqgnh$-modules
\cite{Jim94}:
\bea
\Phi^{\mu V}_\l(z)&:&~~ V(\l)\longrightarrow V(\mu)\otimes V_z,
     \label{Phi}\\
\Phi^{\mu V^*}_\l(z)&:&~~ V(\l)\longrightarrow V(\mu)\otimes V^{*S}_z,
     \label{Phi*}\\
\Psi^{V\mu}_\l(z)&:&~~ V(\l)\longrightarrow V_z\otimes V(\mu),
     \label{Psi}\\
\Psi^{V^*\mu}_\l(z)&:&~~ V(\l)\longrightarrow V^{*S}_z\otimes V(\mu).
     \label{Psi*}
\eea
They are intertwiners in the sense that for any $x\in \uqgnh$
\beq
\Xi(z)\cdot x=\D(x)\cdot\Xi(z),~~~~\Xi(z)=\Phi^{\mu V}_\l(z),~
   \Phi^{\mu V^*}_\l(z),~\Psi^{V\mu}_\l(z),~\Psi^{V^*\mu}_\l(z).
   \label{intertwiner1}
\eeq
These intertwiners are even operators, that is their gradings are
$[\Phi^{\mu V}_\l(z)] = [\Phi^{\mu V^*}_\l(z)] = [\Psi^{V\mu}_\l(z)]
= [\Psi^{V^*\mu}_\l(z)] = 0$. According to \cite{Jim94},
$\Phi^{\mu V}_\l(z)~\lt(\Phi^{\mu V^*}_\l(z)\rt)$ is called type I
(dual) vertex operator and
$\Psi^{V\mu}_\l(z)~\lt(\Psi^{V^*\mu}_\l(z)\rt)$ type II (dual) vertex
operator.
The vertex operators can be expanded in the form 
\bea
&&\Phi^{\mu V}_\l(z)=\sum_{j=1}^{2N}\,\Phi^{\mu V}_{\l,j}(z)
   \otimes v_j, ~~~~
\Phi^{\mu V^*}_\l(z)=\sum_{j=1}^{2N}\,\Phi^{\mu V^*}_{\l,j}(z)
   \otimes v^*_j,\\ 
&&\Psi^{V\mu}_\l(z)=\sum_{j=1}^{2N}\,v_j\otimes \Psi^{V\mu}_{\l,j}(z),~~~~
\Psi^{V^*\mu}_\l(z)=\sum_{j=1}^{2N}\,v^*_j\otimes \Psi^{V^*\mu}_{\l,j}(z).
\eea

Following \cite{Zha98}, we  introduce the combinations of 
bosonic oscillators:
\bea
&&A^{*i}_n=\sum_{l=1}^{2N-1}a^{-1}_{il}A^l_n+\frac{2}{q^n+q^{-n}}  
  a^{-1}_{i,2N}A^{2N}_n,~~~i=1,2,\cdots,2N-1,\\
& &A^{*2N}_n=2N\sum^{2N}_{l=1}a^{-1}_{2N,l}A^{l}_n,~~~ n\neq 0,\\
&&A^{*j}_0=\sum_{l=1}^{2N}a^{-1}_{jl}A^l_0,~~~~
  Q^*_{A^j}=\sum_{l=1}^{2N}a^{-1}_{jl}Q_{A^l},~~~~j=1,2,\cdots,2N,
\eea
\noindent and define the currents
\bea
&&H^{*,j}(z;\k)=Q^*_{A^j}+A^{*j}_0\ln z-\sum_{n\neq
0}\frac{A^{*j}_n}{[n]_q} q^{k|n|}z^{-n},~~~j=1,2,\cdots ,2N-1,\\
&&B_{2N}(z;\k)=Q^*_{A^{2N}}+A^{*2N}_0\ln z-\sum_{n\neq 0}
\frac{A^{*2N}_n}{[n]_q}q^{\k |n|}z^{-n},\\
&&B_{1}(z;\k)=Q^*_{A^{2N}}+A^{*2N}_0\ln z+\frac{N-1}{N}\sum_{n\neq 0}
\frac{A^{*2N}_n}{[n]_q}q^{\k |n|}z^{-n}.
\eea
Consider  the even  operators 
\bea
&&\phi(z)=\sum_{j=1}^{2N}\,\phi_j(z)
   \otimes v_j, ~~~~
\phi^*(z)=\sum_{j=1}^{2N}\,\phi^*_j(z)
   \otimes v^*_j,\no\\
&&\psi(z)=\sum_{j=1}^{2N}\,v_j\otimes \psi_j(z),~~~~
\psi^*(z)=\sum_{j=1}^{2N}\,v^*_j\otimes \psi^*_j(z).\no
\eea
They obey the same intertwining relations as $\Phi_\l^{\mu V}(z),~
\Phi_\l^{\mu V^*}(z),~ \Psi_\l^{V\mu}(z)$ and $\Psi_\l^{V^*\mu}(z)$,
respectively.

Intertwining operators which satisfy (\ref{intertwiner1}) for any 
$x\in\uqsnh$ have been constructed by one of the authors \cite{Zha98}.
We extend the result to $\uqgnh$ by requiring that the vertex 
operators found in \cite{Zha98}
also satisfy (\ref{intertwiner1}) for the element 
$x=H^{2N}_m$ (which extends $\uqsnh$ to $\uqgnh$). We find


\bea
&&\phi_{2N}(z)=:e^{-H^{*2N-1}(qz;\frac{1}{2})+B_{2N}(zq^2;\frac{1}{2})}
e^{c^{N}(qz)}:e^{-\sqrt{-1}\pi N_f},\no\\
&&(-1)^l\phi_l(z)=[\phi_{l+1}(z)~,~f_l]_{q^{(-1)^l}},\no\\
&&\phi^*_1(z)=:e^{H^{*,1}(qz;\frac{1}{2})+B_1(qz;\frac{1}{2})}:
e^{\sqrt{-1}\pi N_f},\no\\
&&q^{(-1)^{l+1}}\phi^*_{l+1}(z)=[\phi^*_l(z)~,f_l]_{q^{(-1)^{l+1}}},\no\\
&&\psi_1(z)=:e^{-H^{*,1}(qz;-\frac{1}{2})-B_1(qz;-\frac{1}{2})}:
e^{-\sqrt{-1}\pi N_f},\no\\
&&\psi_{l+1}(z)=[\psi_l(z)~,e_l]_{q^{(-1)^{l+1}}},\no\\
&&\psi^*_{2N}(z)=:e^{H^{*2N-1}(qz;-\frac{1}{2})-B_{2N}(z;-\frac{1}{2})}
\partial_z\{e^{-c^{N}(qz)}\}:e^{\sqrt{-1}\pi N_f},\no\\            
&&(-1)^{l+1}q^{(-1)^l}\psi^*_l(z)=[\psi^*_{l+1}(z)~,~e_l]_{q^{(-1)^l}}.
\label{Vertex-operator}
\eea
where the $q$-difference derivative $\partial_z$ is defined by
\bea
\partial_zf(z)=\frac{f(zq)-f(zq^{-1})}{(q-q^{-1})z}\no
\eea
and $N_f=\sum_{l=1}^{N}a^{2l}_0$
 is the Fermi-number operator satisfying
\begin{eqnarray*}
(-1)^{N_f}\Xi(z)=(-1)^{[\Xi(z)]}\Xi(z),~~~~{\rm for}~ \Xi(z)=X^{\pm,i}(z), 
\;\phi_i(z),\;\phi^*_i(z),\;\psi_i(z),\;\psi^*_i(z).
\end{eqnarray*}
Since $\Phi^{\mu V}_{\l}(z)$, $\Phi^{\mu V^*}_{\l}(z)$, 
$\Psi^{V\mu}_{\l}(z)$ and $\Psi^{V^*\mu}_{\l}(z)$ obey the same
intertwining relations as $\phi(z)$, $\phi^*(z)$, $\psi(z)$ and 
$\psi^*(z)$ respectively, we have 
\begin{Proposition}:
The vertex operators $\Phi^{\mu V}_\l(z),\;
\Phi^{\mu V^*}_\l(z),\; \Psi^{V\mu}_\l(z)$ and $\Psi^{V^*\mu}_\l(z)$, if
they exist,
have  the same bosonized expressions as the operators
$\phi(z),\;\phi^*(z),\; \psi(z)$
and $\psi^*(z)$, respectively. 
\end{Proposition}

\section{Exchange Relations of Vertex Operators}
In this section, we derive the exchange relations of the type I and type
II bosonized vertex operators of $\uqgnh$. As expected, these 
vertex operators satisfy the graded Faddeev-Zamolodchikov algebra.

\subsection{The R-matrix}
Let $R(z) \in End(V\otimes V)$ be the R-matrix of $\uqgnh$, defined by
:
\bea
R(z)(v_i\otimes v_j)=\sum_{k,l=1}^{2N}R^{ij}_{kl}(z)v_k\otimes v_l,~~
\forall v_i,v_j,v_k,v_l\in V,
\eea
\noindent where 
\begin{eqnarray*}
& &R^{2l-1,2l-1}_{2l-1,2l-1}(z)=1,~~~~
R^{2l,2l}_{2l,2l}(z)=\frac{zq^{-1}-q}{zq-q^{-1}},~~~ l=1,2,\cdots,N,\\
& &R^{ij}_{ij}(z)=\frac{z-1}{zq-q^{-1}},~~~i\neq j,\\
& &R^{ji}_{ij}(z)=\frac{q-q^{-1}}{zq-q^{-1}}(-1)^{[i][j]},~~~i<j,\\
& &R^{ji}_{ij}(z)=\frac{(q-q^{-1})z}{zq-q^{-1}}(-1)^{[i][j]},~~~i>j,\\
& &R^{ij}_{kl}(z)=0,~~~{\rm otherwise}.
\end{eqnarray*}
\noindent  The R-matrix satisfies the graded Yang-Baxter equation  on 
 $V\otimes V\otimes V$
\begin{eqnarray*}
R_{12}(z)R_{13}(zw)R_{23}(w)=R_{23}(w)R_{13}(zw)R_{12}(z),
\end{eqnarray*}
\noindent and moreover enjoys : (i) initial condition, $R(1)=P$ with 
$P$ being the graded permutation operator; (ii) unitarity condition, 
$R_{12}(\frac{z}{w})R_{21}(\frac{w}{z})=1$ , where $R_{21}(z)=P
R_{12}(z)P$; and (iii) crossing-unitarity,
\begin{eqnarray*}
R^{-1,st_1}(z)R(z)^{st_1}=\frac{(z-1)^2}{(q^{-1}z-q)(zq-q^{-1})}.
\end{eqnarray*}
\noindent The various supertranspositions of the R-matrix are given by 
\begin{eqnarray*}
& &(R^{st_1}(z))^{kl}_{ij}=R^{il}_{kj}(z)(-1)^{[i]([i]+[k])},~~~~
(R^{st_2}(z))^{kl}_{ij}=R^{kj}_{il}(z)(-1)^{[j]([l]+[j])},\\
& &(R^{st_{12}}(z))^{kl}_{ij}=R^{ij}_{kl}(z)
(-1)^{([i]+[j])([i]+[j]+[k]+[l])}=R^{ij}_{kl}(z).
\end{eqnarray*}

\subsection{The graded Faddeev-Zamolodchikov algebra}
We calculate the exchange relations of the type I and type 
II bosonic vertex operators of $\uqgnh$ given by  (\ref{Vertex-operator}).
Define 
\begin{eqnarray*}
\oint dzf(z)=Res(f)=f_{-1},\;{\rm for~ a~ formal~ series~ function} 
~ f(z)=\sum_{n\in \Z}f_nz^n.
\end{eqnarray*}
\noindent Then, the Chevalley generators of $\uqgnh$ can be expressed by
the integrals
\begin{eqnarray*}
e_i=\oint dz X^{+,i}(z),~~~~f_i=\oint dz X^{-,i}(z),~~~i=1,2,\cdots,
2N-1.
\end{eqnarray*}
\noindent One can also get the integral expressions of the bosonic vertex
operators $\phi(z)$, $\phi^*(z)$, $\psi(z)$ and $\phi^*(z)$ from  
(\ref{Vertex-operator}). Using these integral expressions and the
relations  given in appendices A and B, we arrive at  
\begin{Proposition}
: The bosonic vertex operators defined in (\ref{Vertex-operator})
satisfy the graded Faddeev-Zamolodchikov algebra
\bea
&&\phi_j(z_2)\phi_i(z_1)=(\frac{z_2}{z_1})^{2-\frac{1}{N}}
\sum_{k,l=1}^{2N}R^{kl}_{ij}(\frac{z_1}{z_2})
\phi_k(z_1)\phi_l(z_2)(-1)^{[i][j]},\label{ZF1}\\
&&\psi^*_i(z_1)\psi^*_j(z_2)=(\frac{z_1}{z_2})^{2-\frac{1}{N}}
\sum_{k,l=1}^{2N}R^{ij}_{kl}(\frac{z_1}{z_2})
\psi^*_l(z_2)\psi^*_k(z_1)(-1)^{[i][j]},\\
&&\psi^*_i(z_1)\phi_j(z_2)=(\frac{qz_2}{z_1})^{2-\frac{1}{N}}
\phi_j(z_2)\psi^*_i(z_1)(-1)^{[i][j]}.\label{ZF2}
\eea
\end{Proposition}

In the derivation of this proposition the fact that $R^{kl}_{ij}(z)
(-1)^{[k][l]}=R^{kl}_{ij}(z)(-1)^{[i][j]}$ is helpful.

We can also generalize the Miki's construction to the  $\uqgnh$ case. Define 
\begin{eqnarray*}
& &L^+(z)^j_i=\phi_i(zq^{1/ 2})\psi^*_j(zq^{-1/2}),\\
& &L^-(z)^j_i=\phi_i(zq^{-1/ 2})\psi^*_j(zq^{1/2}).
\end{eqnarray*}
\begin{Proposition}
: The L-operators $L^{\pm}(z)$ defined above give a realization of the
super RS algebra \cite{Zha97} at level one  for 
the quantum affine superalgebra $\uqgnh$ 
\begin{eqnarray*}
& &R(z/w)L^{\pm}_1(z)L^{\pm}_2(w)=L^{\pm}_2(w)L^{\pm}_1(z)R(z/w),\\
& &R(z^{+}/w^{-})L^{+}_1(z)L^{-}_2(w)=L^{-}_2(w)L^{+}_1(z)R(z^-/w^+),
\end{eqnarray*}
\noindent where $L^{\pm}_1(z)=L^{\pm}(z)\otimes 1$,
$L^{\pm}_2(z)=1\otimes L^{\pm}(z)$ and $z^{\pm}=zq^{\pm\frac{1}{2}}$.
\end{Proposition}

\noindent {\it Proof}. Straightforward computation
by using (\ref{ZF1})-(\ref{ZF2}).

\sect{Highest Weight Representation of  $\uqg2h$ and 
Associated Vertex Operators}

In this section we study in details the irreducible highest weight 
representations of $\uqg2h$ and their associated vertex operators.

\subsection{Highest weight $\uqg2h$-modules}
We begin by defining the Fock module. Denote by
$F_{\l_1,\l_2,\l_3,\l_4;\l_5,\l_6}$ the bosonic Fock
spaces generated by $a_{-m}^i,c^l_{-m} (m>0)$
 over the vector $|\l_1,\l_2,\l_3,\l_4;\l_5,\l_6>$:
\begin{eqnarray*}
F_{\l_1,\l_2,\l_3,\l_4;\l_5,\l_6}&=&
{\bf C}\lt[a^1_{-1},a^1_{-2}\cdots;a^{2}_{-1},a^{2}_{-2}\cdots; 
a^3_{-1},a^3_{-2},\cdots;a^4_{-1},a^4_{-2},\cdots;\rt.\\
& &\lt.c^1_{-1},c^1_{-2},\cdots;c^{2}_{-1},c^{2}_{-2},\cdots\rt]
|\l_1,\l_2,\l_3,\l_4;\l_5,\l_6>,
\end{eqnarray*}
\noindent where 
\begin{eqnarray*}
|\l_1,\l_2,\l_3,\l_4;\l_5,\l_6>=
e^{\sum_{i=1}^{4}\l_iQ_{a^i}+\l_5Q_{c^1}+\l_6Q_{c^2}}|0>.
\end{eqnarray*}
\noindent The vacuum vector $|0>$ is defined by 
$a^i_m|0>=c^l_m|0>=0$ for $i=1,2,3,4$, $\;l=1,2$ and $m\geq 0$. Obviously, 
\begin{eqnarray*}
&& a^i_m|\l_1,\l_2,\l_3,\l_4;\l_5,\l_6>=0,~~~{\rm for}~i=1,2,3,4~{\rm
   and}~ m>0,\\
&& c^l_m|\l_1,\l_2,\l_3,\l_4;\l_5,\l_6>=0,~~~{\rm for}~l=1,2~{\rm and}~
   m>0.
\end{eqnarray*}
To obtain the highest weight vectors of $\uqg2h$, we impose the
conditions:
\bea
&&e_i|\l_1,\l_2,\l_3,\l_4;\l_5,\l_6>=0,~~~ i=0,1,2,3,\no\\
&&h_j|\l_1,\l_2,\l_3,\l_4;\l_5,\l_6>=
  \l^j|\l_1,\l_2,\l_3,\l_4;\l_5,\l_6>,~~~
 j=0,1,2,3,4.
\eea
\noindent Solving these equations, we obtain the following classification:
\begin{enumerate}
\item $(\l_1,\l_2,\l_3,\l_4;\l_5,\l_6)=(\b,-\b,\b,-\b;0,0)$, where $\b$ is
arbitrary. The weight of this vector is $(\l^0,\l^1,\l^2,\l^3,\l^4)$=$
(1,0,0,0,4\b)$. We have $|\L_0+4\b\L_4>$=$|\b,-\b,\b,-\b;0,0>$.

\item $(\l_1,\l_2,\l_3,\l_4;\l_5,\l_6)=(\b+1,-\b-1,\b+1,-\b;0,0)$, where
$\b$ is arbitrary. The weight of this vector is 
$(\l^0,\l^1,\l^2,\l^3,\l^4)$=$
(0,0,0,1,4\b+3)$. We have
$|\L_3+(4\b+3)\L_4>$=$|\b+1,-\b-1,\b+1,-\b;0,0>$.

\item $(\l_1,\l_2,\l_3,\l_4;\l_5,\l_6)=(\b+1,-\b-1+\a,\b,-\b;-\a,0)$,
where
$\b$ is arbitrary. The weight of this vector is 
$(\l^0,\l^1,\l^2,\l^3,\l^4)$=$
(0,\a,1-\a,0,4\b+2-\a)$. We have
$|\a\L_1+(1-\a)\L_2+(4\b+2-\a)\L_4>$=$|\b+1,-\b-1+\a,\b,-\b;-\a,0>$.
\end{enumerate}
According to this classification, we introduce the Fock  
spaces
\bea 
&&\F_{((0,1);\b)}=\oplus_{i,j,k\in
\Z}F_{\b+i,\; -\b-i+j,\; \b-j+k,\; -\b-k;\;
i-j,\; k},\label{Fack-space1}\\
&&\F_{((1,0);\b)}=\oplus_{i,j,k\in \Z}
F_{\b+1+i,\; -\b-1-i+j,\; \b+1-j+k,\; -\b-k;\; i-j,\; k},\\
&&\F_{(\a;\b)}=\oplus_{i,j,k\in \Z}
F_{\b+1+i,\; -\b-1+\a-i+j,\;
\b-j+k,\; -\b-k;\; -\a+i-j,\; k}.\label{Fack-space2}
\eea
It can be shown that the bosonized action of $\uqg2h$ on $\F_{(*;\b)}$ 
 is closed, i.e.
\begin{eqnarray*}
&&\uqg2h \F_{(*;\b)}=\F_{(*;\b)},~~~
 {\rm for}~ *=(0,1),(1,0),\a.
\end{eqnarray*}
\noindent Hence each Fack space in 
(\ref{Fack-space1})-(\ref{Fack-space2}) 
constitutes a $\uqg2h$-module. However, these modules are  not 
irreducible in general. To obtain the irreducible representations, we 
introduce two pairs of ``fermionic" currents 
\begin{eqnarray*}
\eta^i(z)=\sum_{n\in \Z}\eta_n^iz^{-n-1}=:e^{c^i(z)}:,~~~~
\xi^i(z)=\sum_{n\in \Z}\xi_n^iz^{-n}=:e^{-c^i(z)}:,~~~i=1,2.
\end{eqnarray*}
\noindent The mode expansion of $\eta^i(z)$ and $\xi^i(z)$ is well 
defined on $\F_{((0,1);\b)}$,$\F_{((1,0);\b)}$ and $\F_{(\a;\b)}$ with
$\a\in \Z$, and 
the modes satisfy the relations
\bea
&&\xi^i_m\xi^i_n+\xi^i_n\xi^i_m=\eta^i_m\eta^i_n+\eta^i_n\eta^i_m=0~,~~
\xi^i_m\eta^i_n+\eta^i_n\xi^i_m=\delta_{m+n,0}~,~~i=1,2,\\
&&\xi_m^1\xi^2_n-\xi^2_n\xi^1_m=\xi_m^1\eta^2_n-\eta^2_n\xi^1_m=
\eta_m^1\xi^2_n-\xi^2_n\eta^1_m=\eta_m^1\eta^2_n-\eta^2_n\eta^1_m=0.
\eea
\noindent Thus, we have the direct sum decompositions
\bea
\F_{(*;\b)}=\eta^1_0\xi^1_0\eta^2_0\xi^2_0\F_{(*;\b)}\oplus
\eta^1_0\xi^1_0\xi^2_0\eta^2_0\F_{(*;\b)}\oplus
\xi^1_0\eta^1_0\eta^2_0\xi^2_0\F_{(*;\b)}\oplus
\xi^1_0\eta^1_0\xi^2_0\eta^2_0\F_{(*;\b)},
\eea
\noindent with $*=(0,1),(1,0),\a\in \Z$. As usual, we name 
\begin{eqnarray*}
KKer|_{\F_{(*;\b)}}~&{\rm as}&~~
\eta^1_0\xi^1_0\eta^2_0\xi^2_0\F_{(*;\b)},\\
KCoker|_{\F_{(*;\b)}}~&{\rm as}&~~ \eta^1_0\xi^1_0\xi^2_0\eta^2_0
\F_{(*;\b)},\\
CKer|_{\F_{(*;\b)}}~&{\rm as}&~~ \xi^1_0\eta^1_0\eta^2_0\xi^2_0
\F_{(*;\b)},\\
CCoker|_{\F_{(*;\b)}}~&{\rm as}&~~ \xi^1_0\eta^1_0\xi^2_0\eta^2_0
\F_{(*;\b)}.
\end{eqnarray*}
\noindent  Since $\eta^1_0$ and $\eta^2_0$ commute (or
anticommute) with the bosonized 
actions of $\uqg2h$, $KKer|_{\F_{(*;\b)}}$, $CKer|_{\F_{(*;\b)}}$,
 $KCoker|_{\F_{(*;\b)}}$
 and $CCoker|_{\F_{(*;\b)}}$ are all the $\uqg2h$-modules.

{}From now on, we study the characters and supercharacters  of 
these $\uqg2h$-modules which are constructed in the bosonic Fock spaces.
We first of all bosonize the derivation operator $d$ as 
\bea
d&=&-\sum_{m>0}\frac{m^2}{[m]_q^2}\{
\sum_{i=1}^{3}A^i_{-m}A^{*i}_{m}+\frac{1}{2(q^m+q^{-m})}A^4_{-m}A^{*4}_{m}
+c^1_{-m}c^1_m+c^2_{-m}c^2_m\}\no\\
&&~~-\frac{1}{2}\{\sum_{i=1}^{4}A^i_0A^{*i}_0+c^1_0(c^1_0+1)+c^2_0(c^2_0+1)\}.
\label{d}
\eea
\noindent One can easily check that this $d$ obeys the commutation
relations 
\begin{eqnarray*}
[d,h_j]=0,~~~ [d,h^j_m]=mh^j_m,~~~ [d,X^{\pm,i}_m]=mX^{\pm,i}_m, ~~
j=1,2,3,4,~~ i=1,2,3,
\end{eqnarray*}
\noindent as required. Moreover, we have $[d,\xi^l_0]=[d,\eta^l_0]=0$ for 
$l=1,2.$

The character and supercharacter of a $\uqg2h$-module $M$ are defined by 
\bea
Ch_{M}(q,x_1,x_2,x_3,x_4)&=&tr_M(q^{-d}x_1^{h_1}x_2^{h_2}x_3^{h_3}
  x_4^{h_4}),\\
Sch_{M}(q,x_1,x_2,x_3,x_4)&=&Str_M(q^{-d}x_1^{h_1}x_2^{h_2}x_3^{h_3}
  x_4^{h_4})\no\\
&=&tr_M((-1)^{N_f}q^{-d}x_1^{h_1}x_2^{h_2}x_3^{h_3}
  x_4^{h_4}),
\eea
\noindent respectively.

\begin{itemize}
\item  (I) {\it Character of $\F_{(\a;\b)}$ for $\a\not\in \Z$}. Since  
$\eta^1_0,\eta^2_0$ are  not defined on this module, it is expected 
that $\F_{(\a;\b)}$ is an  irreducible highest weight $\uqg2h$-module.
Thus, we have
\end{itemize}
\vspace{0.4truecm}
\noindent {\bf Conjecture 1} : {\it We have the identification 
of the highest weight $\uqg2h$-modules:
\bea
\F_{(\a;\b)}\cong V(\a\L_1+(1-\a)\L_2+(4\b+2-\a)\L_4)
~~for ~~\a\not\in \Z~ and~arbitrary ~\b~,\no
\eea
where and throughout $V(\l)$ denotes the irreducible highest
weight
$\uqg2h$-module with the highest weight $\l$.}

\begin{Proposition}
: The character and supercharacter of $\F_{(\a;\b)}$ are  
\bea
&&Ch_{\F_{(\a;\b)}}(q,x_1,x_2,x_3,x_4)=\frac{q^{\frac{1}{2}\a(2\b+1)}}
{\prod_{n=1}^{\infty}(1-q^n)^6}~~~~~~~~~~~~~~~~\no\\
&&~~~~~~~~~~\times \sum_{i,j,k\in
\Z}q^{\frac{1}{2}(i^2+j^2+k^2-2kj+i+j+k)}
x_1^{\a+j}x_2^{1-\a+i-k}x_3^{-j}x_4^{(4\b+2-\a+2i-2j+2k)},\no\\
&&Sch_{\F_{(\a;\b)}}(q,x_1,x_2,x_3,x_4)=\frac{q^{\frac{1}{2}\a(2\b+1)}}
{\prod_{n=1}^{\infty}(1-q^n)^6}~~~~~~~~~~~~~~~~\no\\
&&~~~~~~~~~~\times \sum_{i,j,k\in
\Z}(-1)^{\a-1+i-j+k}q^{\frac{1}{2}(i^2+j^2+k^2-2kj+i+j+k)}
x_1^{\a+j}x_2^{1-\a+i-k}x_3^{-j}x_4^{(4\b+2-\a+2i-2j+2k)}.\no
\eea
\end{Proposition}
\begin{itemize}
\item (II) {\it Characters and supercharacters of $KKer_{\F_{(*;\b)}}$,
 $CKer_{\F_{(*;\b)}}$, $KCoker_{\F_{(*;\b)}}$
and $CCoker_{\F_{(*;\b)}}$ for $ *=(0,1),(1,0),\a\in \Z$}.
In this case, $\eta^1_0$ and $\eta^2_0$ are
 well defined. We shall calculate the  characters and
supercharacters  of these modules by using the BRST
resolution \cite{Yan99}.
\end{itemize}
Let us define the Fock spaces, for $l_1,l_2\in \Z$
\begin{eqnarray*}
&&\F^{(l_1,l_2)}_{((0,1);\b)}=\oplus_{i,j,k\in \Z}F_{\b+i,\; -\b-i+j,
\; \b-j+k,\; -\b-k;\; i-j+l_1,\; k+l_2},\\
&&\F^{(l_1,l_2)}_{((1,0);\b)}=\oplus_{i,j,k\in \Z}F_{\b+1+i,\;
-\b-1-i+j,\; 
\b+1-j+k,\; -\b-k;\; i-j+l_1,\; k+l_2},\\
&&\F^{(l_1,l_2)}_{(\a;\b)}=\oplus_{i,j,k\in \Z}F_{\b+1+i,\;
-\b-1+\a-i+j,\; 
\b-j+k,\; -\b-k;\; -\a+i-j+l_1,\; k+l_2}.
\end{eqnarray*}
\noindent We have $\F^{(0,0)}_{(*;\b)}=\F_{(*;\b)}$. It can be shown that  
$\eta^i_0$ and $\xi^i_0$ intertwine these Fock spaces in the following
fashions
\begin{eqnarray*}
&&\eta^1_0:~~\F^{(l_1,l_2)}_{(*;\b)}\longrightarrow 
\F^{(l_1+1,l_2)}_{(*;\b)},~~~~
\eta^2_0:~~\F^{(l_1,l_2)}_{(*;\b)}\longrightarrow 
\F^{(l_1,l_2+1)}_{(*;\b)},\\
&&\xi^1_0:~~\F^{(l_1,l_2)}_{(*;\b)}\longrightarrow 
\F^{(l_1-1,l_2)}_{(*;\b)},~~~~
\xi^2_0:~~\F^{(l_1,l_2)}_{(*;\b)}\longrightarrow
\F^{(l_1,l_2-1)}_{(*;\b)}.
\end{eqnarray*}
\noindent We have the following two BRST complexes:
\bea
\begin{array}{ccccccc}
\cdots &\stackrel{Q^{(1)}_{l_1-1}=\eta^1_0}{\longrightarrow}&
\F^{(l_1,l_2)}_{(*;\b)}&\stackrel{Q^{(1)}_{l_1}=\eta^1_0}{\longrightarrow}&
\F^{(l_1+1,l_2)}_{(*;\b)}&\stackrel{Q^{(1)}_{l_1+1}=\eta^1_0}{\longrightarrow}&
\cdots\\
&&|{\bf O}&&|{\bf O}&&\\
\cdots &\stackrel{Q^{(1)}_{l_1-1}=\eta^1_0}{\longrightarrow}&
\F^{(l_1,l_2)}_{(*;\b)}&\stackrel{Q^{(1)}_{l_1}=\eta^1_0}{\longrightarrow}&
\F^{(l_1+1,l_2)}_{(*;\b)}&\stackrel{Q^{(1)}_{l_1+1}=\eta^1_0}{\longrightarrow}&
\cdots
\end{array}
\eea
\noindent and 
\bea
\begin{array}{ccccccc}
\cdots &\stackrel{Q^{(2)}_{l_2-1}=\eta^2_0}{\longrightarrow}&
\F^{(l_1,l_2)}_{(*;\b)}&\stackrel{Q^{(2)}_{l_2}=\eta^2_0}{\longrightarrow}&
\F^{(l_1,l_2+1)}_{(*;\b)}&\stackrel{Q^{(2)}_{l_2+1}=\eta^2_0}{\longrightarrow}&
\cdots\\
&&|{\bf O}&&|{\bf O}&&\\
\cdots &\stackrel{Q^{(2)}_{l_2-1}=\eta^2_0}{\longrightarrow}&
\F^{(l_1,l_2)}_{(*;\b)}&\stackrel{Q^{(2)}_{l_2}=\eta^2_0}{\longrightarrow}&
\F^{(l_1,l_2+1)}_{(*;\b)}&\stackrel{Q^{(2)}_{l_2+1}=\eta^2_0}{\longrightarrow}
&
\cdots
\end{array}
\eea
\noindent where ${\bf O}$ is an operator such that 
$\F^{(l_1,l_2)}_{(*;\b)}\longrightarrow \F^{(l_1,l_2)}_{(*;\b)}$, and
moreover 
${\bf O}$ commutes with the BRST charges $Q^{(1)}_l$ and $Q^{(2)}_{l}$. 
Then,   
\begin{Proposition}:\label{brst}
\bea
&&Ker_{Q^{(i)}_l}=Im_{Q^{(i)}_{l-1}},~~i=1,2, ~ ~{\rm for~ any~} l\in \Z,
\no\\
&&tr({\bf O})|_{Ker_{Q^{(i)}_l}}=tr({\bf O})|_{Im_{Q^{(i)}_{l-1}}}
=tr({\bf O})|_{Coker_{Q^{(i)}_{l-1}}}~~.
\eea
\end{Proposition}
\noindent {\it Proof}. It follows from the fact that 
$\eta^i_0\xi^i_0+\xi^i_0\eta^i_0=1$, $(\eta^i_0)^2=(\xi^i_0)^2=0$ and 
$\eta^i_0\xi^i_0$ ($\xi^i_0\eta^i_0$) are  the projection operators  
from $\F^{(l_1,l_2)}_{(*;\b)}$ to $Ker_{Q^{(i)}_{l_i}}$ 
($Coker_{Q^{(i)}_{l_i}}$).

\vspace{0.5truecm}

By proposition \ref{brst}, we can compute   
the characters and 
supercharacters of $KKer_{\F_{(*;\b)}}$, $CKer_{\F_{(*;\b)}}$, 
 $KCoker_{\F_{(*;\b)}}$ and $CCoker_{\F_{(*;\b)}}$ for 
$ *=(0,1),(1,0), \a\in \Z$. We have 
\begin{Proposition}
: The characters and
supercharacters of $KKer_{\F_{(0,1);\b)}}$, $CKer_{\F_{((0,1);\b)}}$,\\
$KCoker_{\F_{((0,1);\b)}}$ and $CCoker_{\F_{((0,1);\b)}}$  are given by
\bea
&&Ch_{KKer_{\F_{((0,1);\b)}}}(q,x_1,x_2,x_3,x_4)
=\frac{1}{\prod_{n=1}^{\infty}(1-q^n)^6}   
\sum_{l_1,l_2=1}^{\infty}(-1)^{l_1+l_2}q^{\frac{1}{2}(l_1^2+l_2^2-
l_1-l_2)}\no\\
&&~~~~~~~~~~~~~~~~\times \sum_{i,j,k\in
\Z}q^{\frac{1}{2}(i^2+j^2+k^2-2kj+(1-2l_1)i-(1-2l_1)j+(1-2l_2)k)}
x_1^{j}x_2^{i-k}x_3^{-j}x_4^{(4\b+2i-2j+2k)},\no\\
&&Ch_{CKer_{\F_{((0,1);\b)}}}(q,x_1,x_2,x_3,x_4)
=\frac{1}{\prod_{n=1}^{\infty}(1-q^n)^6}   
\sum_{l_1,l_2=1}^{\infty}(-1)^{l_1+l_2}q^{\frac{1}{2}(l_1^2+l_2^2+
l_1-l_2)}\no\\
&&~~~~~~~~~~~~~~~~\times \sum_{i,j,k\in
\Z}q^{\frac{1}{2}(i^2+j^2+k^2-2kj+(1+2l_1)i-(1+2l_1)j+(1-2l_2)k)}
x_1^{j}x_2^{i-k}x_3^{-j}x_4^{(4\b+2i-2j+2k)},\no\\
&&Ch_{KCoker_{\F_{((0,1);\b)}}}(q,x_1,x_2,x_3,x_4)
=\frac{1}{\prod_{n=1}^{\infty}(1-q^n)^6}   
\sum_{l_1,l_2=1}^{\infty}(-1)^{l_1+l_2}q^{\frac{1}{2}(l_1^2+l_2^2-
l_1+l_2)}\no\\
&&~~~~~~~~~~~~~~~~\times \sum_{i,j,k\in
\Z}q^{\frac{1}{2}(i^2+j^2+k^2-2kj+(1-2l_1)i-(1-2l_1)j+(1+2l_2)k)}
x_1^{j}x_2^{i-k}x_3^{-j}x_4^{(4\b+2i-2j+2k)},\no\\
&&Ch_{CCoker_{\F_{((0,1);\b)}}}(q,x_1,x_2,x_3,x_4)
=\frac{1}{\prod_{n=1}^{\infty}(1-q^n)^6}   
\sum_{l_1,l_2=1}^{\infty}(-1)^{l_1+l_2}q^{\frac{1}{2}(l_1^2+l_2^2+
l_1+l_2)}\no\\
&&~~~~~~~~~~~~~~~~\times \sum_{i,j,k\in
\Z}q^{\frac{1}{2}(i^2+j^2+k^2-2kj+(1+2l_1)i-(1+2l_1)j+(1+2l_2)k)}
x_1^{j}x_2^{i-k}x_3^{-j}x_4^{(4\b+2i-2j+2k)},\no
\eea
\noindent and 
\bea
&&Sch_{KKer_{\F_{((0,1);\b)}}}(q,x_1,x_2,x_3,x_4)
=\frac{1}{\prod_{n=1}^{\infty}(1-q^n)^6}
\sum_{l_1,l_2=1}^{\infty}(-1)^{l_1+l_2}q^{\frac{1}{2}(l_1^2+l_2^2-
l_1-l_2)}\no\\
&&~~~~~~~~~~~~~~~~\times \sum_{i,j,k\in \Z}\lt[
(-1)^{i-j+k}q^{\frac{1}{2}(i^2+j^2+k^2-2kj
+(1-2l_1)i-(1-2l_1)j+(1-2l_2)k)}\rt.\no\\
&&~~~~~~~~~~~~~~~~\times
\lt.x_1^{j}x_2^{i-k}x_3^{-j}x_4^{(4\b+2i-2j+2k)}\rt],\no\\
&&Sch_{CKer_{\F_{((0,1);\b)}}}(q,x_1,x_2,x_3,x_4)
=\frac{1}{\prod_{n=1}^{\infty}(1-q^n)^6}
\sum_{l_1,l_2=1}^{\infty}(-1)^{l_1+l_2}q^{\frac{1}{2}(l_1^2+l_2^2+
l_1-l_2)}\no\\
&&~~~~~~~~~~~~~~~~\times \sum_{i,j,k\in
\Z}\lt[(-1)^{i-j+k}q^{\frac{1}{2}
(i^2+j^2+k^2-2kj+(1+2l_1)i-(1+2l_1)j+(1-2l_2)k)}\rt.\no\\
&&~~~~~~~~~~~~~~~~~\times 
\lt.x_1^{j}x_2^{i-k}x_3^{-j}x_4^{(4\b+2i-2j+2k)}\rt],\no\\ 
&&Sch_{KCoker_{\F_{((0,1);\b)}}}(q,x_1,x_2,x_3,x_4)
=\frac{1}{\prod_{n=1}^{\infty}(1-q^n)^6}
\sum_{l_1,l_2=1}^{\infty}(-1)^{l_1+l_2}q^{\frac{1}{2}(l_1^2+l_2^2-
l_1+l_2)}\no\\
&&~~~~~~~~~~~~~~~~~\times \sum_{i,j,k\in
\Z}\lt[(-1)^{i-j+k}q^{\frac{1}{2}
(i^2+j^2+k^2-2kj+(1-2l_1)i-(1-2l_1)j+(1+2l_2)k)}\rt.\no\\
&&~~~~~~~~~~~~~~~~~\times 
\lt.x_1^{j}x_2^{i-k}x_3^{-j}x_4^{(4\b+2i-2j+2k)}\rt],\no\\
&&Sch_{CCoker_{\F_{((0,1);\b)}}}(q,x_1,x_2,x_3,x_4)
=\frac{1}{\prod_{n=1}^{\infty}(1-q^n)^6}
\sum_{l_1,l_2=1}^{\infty}(-1)^{l_1+l_2}q^{\frac{1}{2}(l_1^2+l_2^2+
l_1+l_2)}\no\\
&&~~~~~~~~~~~~~~~~~~\times \sum_{i,j,k\in
\Z}\lt[(-1)^{i-j+k}q^{\frac{1}{2}
(i^2+j^2+k^2-2kj+(1+2l_1)i-(1+2l_1)j+(1+2l_2)k)}\rt.\no\\
&&~~~~~~~~~~~~~~~~~~\times 
\lt.x_1^{j}x_2^{i-k}x_3^{-j}x_4^{(4\b+2i-2j+2k)}\rt].\no
\eea
\end{Proposition}

\begin{Proposition}
: The characters and
supercharacters of $KKer_{\F_{(1,0);\b)}}$, $CKer_{\F_{((1,0);\b)}}$,\\
$KCoker_{\F_{((1,0);\b)}}$ and $CCoker_{\F_{((1,0);\b)}}$  are given by 
\bea
&&Ch_{KKer_{\F_{((1,0);\b)}}}(q,x_1,x_2,x_3,x_4)
=\frac{q^{\frac{1}{2}(2\b+1)}}{\prod_{n=1}^{\infty}(1-q^n)^6}
\sum_{l_1,l_2=1}^{\infty}(-1)^{l_1+l_2}q^{\frac{1}{2}(l_1^2+l_2^2-
l_1-l_2)}\no\\
&&~~~~~~~~~~~\times \sum_{i,j,k\in
\Z}q^{\frac{1}{2}(i^2+j^2+k^2-2kj+(1-2l_1)i-(1-2l_1)j+(3-2l_2)k)}
x_1^{j}x_2^{i-k}x_3^{1-j}x_4^{(4\b+3+2i-2j+2k)},\no\\
&&Ch_{CKer_{\F_{((1,0);\b)}}}(q,x_1,x_2,x_3,x_4)
=\frac{q^{\frac{1}{2}(2\b+1)}}{\prod_{n=1}^{\infty}(1-q^n)^6}
\sum_{l_1,l_2=1}^{\infty}(-1)^{l_1+l_2}q^{\frac{1}{2}(l_1^2+l_2^2+   
l_1-l_2)}\no\\
&&~~~~~~~~~~~~~\times \sum_{i,j,k\in
\Z}q^{\frac{1}{2}(i^2+j^2+k^2-2kj+(1+2l_1)i-(1+2l_1)j+(3-2l_2)k)}
x_1^{j}x_2^{i-k}x_3^{1-j}x_4^{(4\b+3+2i-2j+2k)},\no\\
&&Ch_{KCoker_{\F_{((1,0);\b)}}}(q,x_1,x_2,x_3,x_4)
=\frac{q^{\frac{1}{2}(2\b+1)}}{\prod_{n=1}^{\infty}(1-q^n)^6}
\sum_{l_1,l_2=1}^{\infty}(-1)^{l_1+l_2}q^{\frac{1}{2}(l_1^2+l_2^2-
l_1+l_2)}\no\\
&&~~~~~~~~~~~~\times \sum_{i,j,k\in 
\Z}q^{\frac{1}{2}(i^2+j^2+k^2-2kj+(1-2l_1)i-(1-2l_1)j+(3+2l_2)k)}
x_1^{j}x_2^{i-k}x_3^{1-j}x_4^{(4\b+3+2i-2j+2k)},\no\\
&&Ch_{CCoker_{\F_{((1,0);\b)}}}(q,x_1,x_2,x_3,x_4)
=\frac{q^{\frac{1}{2}(2\b+1)}}{\prod_{n=1}^{\infty}(1-q^n)^6}
\sum_{l_1,l_2=1}^{\infty}(-1)^{l_1+l_2}q^{\frac{1}{2}(l_1^2+l_2^2+
l_1+l_2)}\no\\
&&~~~~~~~~~~~~\times \sum_{i,j,k\in
\Z}q^{\frac{1}{2}(i^2+j^2+k^2-2kj+(1+2l_1)i-(1+2l_1)j+(++2l_2)k)}
x_1^{j}x_2^{i-k}x_3^{-j}x_4^{(4\b+3+2i-2j+2k)},\no
\eea
\noindent and 
\bea
&&Sch_{KKer_{\F_{((1,0);\b)}}}(q,x_1,x_2,x_3,x_4)   
=\frac{q^{\frac{1}{2}(2\b+1)}}{\prod_{n=1}^{\infty}(1-q^n)^6}
\sum_{l_1,l_2=1}^{\infty}(-1)^{l_1+l_2}q^{\frac{1}{2}(l_1^2+l_2^2-
l_1-l_2)}\no\\
&&~~~~~~~~~~~~\times \sum_{i,j,k\in
\Z}\lt[(-1)^{1+i-j+k}q^{\frac{1}{2}
(i^2+j^2+k^2-2kj+(1-2l_1)i-(1-2l_1)j+(3-2l_2)k)}\rt.\no\\
&&~~~~~~~~~~~~~\times 
\lt.x_1^{j}x_2^{i-k}x_3^{1-j}x_4^{(4\b+3+2i-2j+2k)}\rt],\no\\
&&Sch_{CKer_{\F_{((1,0);\b)}}}(q,x_1,x_2,x_3,x_4)   
=\frac{q^{\frac{1}{2}(2\b+1)}}{\prod_{n=1}^{\infty}(1-q^n)^6}
\sum_{l_1,l_2=1}^{\infty}(-1)^{l_1+l_2}q^{\frac{1}{2}(l_1^2+l_2^2+
l_1-l_2)}\no\\
&&~~~~~~~~~~~~~\times \sum_{i,j,k\in
\Z}\lt[(-1)^{1+i-j+k}q^{\frac{1}{2}
(i^2+j^2+k^2-2kj+(1+2l_1)i-(1+2l_1)j+(3-2l_2)k)}\rt.\no\\
&&~~~~~~~~~~~~~\times 
\lt.x_1^{j}x_2^{i-k}x_3^{1-j}x_4^{(4\b+3+2i-2j+2k)}\rt],\no\\
&&Sch_{KCoker_{\F_{((1,0);\b)}}}(q,x_1,x_2,x_3,x_4) 
=\frac{q^{\frac{1}{2}(2\b+1)}}{\prod_{n=1}^{\infty}(1-q^n)^6}
\sum_{l_1,l_2=1}^{\infty}(-1)^{l_1+l_2}q^{\frac{1}{2}(l_1^2+l_2^2-
l_1+l_2)}\no\\
&&~~~~~~~~~~~~~\times \sum_{i,j,k\in
\Z}\lt[(-1)^{1+i-j+k}q^{\frac{1}{2}
(i^2+j^2+k^2-2kj+(1-2l_1)i-(1-2l_1)j+(3+2l_2)k)}\rt.\no\\
&&~~~~~~~~~~~~~~\times 
\lt.x_1^{j}x_2^{i-k}x_3^{1-j}x_4^{(4\b+3+2i-2j+2k)}\rt],\no\\
&&Sch_{CCoker_{\F_{((1,0);\b)}}}(q,x_1,x_2,x_3,x_4)
=\frac{q^{\frac{1}{2}(2\b+1)}}{\prod_{n=1}^{\infty}(1-q^n)^6}
\sum_{l_1,l_2=1}^{\infty}(-1)^{l_1+l_2}q^{\frac{1}{2}(l_1^2+l_2^2+
l_1+l_2)}\no\\
&&~~~~~~~~~~~~~~\times \sum_{i,j,k\in
\Z}\lt[(-1)^{1+i-j+k}q^{\frac{1}{2}
(i^2+j^2+k^2-2kj+(1+2l_1)i-(1+2l_1)j+(++2l_2)k)}\rt.\no\\
&&~~~~~~~~~~~~~~\times 
\lt.x_1^{j}x_2^{i-k}x_3^{-j}x_4^{(4\b+3+2i-2j+2k)}\rt].\no
\eea
\end{Proposition}

\begin{Proposition}
: The characters and
supercharacters of $KKer_{\F_{(\a;\b)}}$, $CKer_{\F_{(\a;\b)}}$,\\
$KCoker_{\F_{(\a;\b)}}$ and $CCoker_{\F_{(\a;\b)}}$ 
for $\a\in \Z$ are given by
\bea
&&Ch_{KKer_{\F_{(\a;\b)}}}(q,x_1,x_2,x_3,x_4)
=\frac{q^{\frac{1}{2}\a(2\b+1)}}{\prod_{n=1}^{\infty}(1-q^n)^6}
\sum_{l_1,l_2=1}^{\infty}(-1)^{l_1+l_2}q^{\frac{1}{2}(l_1^2+l_2^2-
(1-2\a)l_1-l_2)}\no\\
&&~~~~~~~~~~~~~~\times \sum_{i,j,k\in
\Z}\lt[q^{\frac{1}{2}(i^2+j^2+k^2-2kj+(1-2l_1)i+(1+2l_1)j+(1-2l_2)k)}\rt.\no\\
&&~~~~~~~~~~~~~~\times
 \lt.x_1^{\a+j}x_2^{1-\a+i-k}x_3^{-j}x_4^{(4\b+2-\a+2i-2j+2k)}\rt],\no\\
&&Ch_{CKer_{\F_{(\a;\b)}}}(q,x_1,x_2,x_3,x_4)
=\frac{q^{\frac{1}{2}\a(2\b+1)}}{\prod_{n=1}^{\infty}(1-q^n)^6}
\sum_{l_1,l_2=1}^{\infty}(-1)^{l_1+l_2}q^{\frac{1}{2}(l_1^2+l_2^2+
(1-2\a)l_1-l_2)}\no\\
&&~~~~~~~~~~~~~~\times \sum_{i,j,k\in
\Z}\lt[q^{\frac{1}{2}(i^2+j^2+k^2-2kj+(1+2l_1)i+(1-2l_1)j+(1-2l_2)k)}\rt.\no\\
&&~~~~~~~~~~~~~~\times 
 \lt.x_1^{\a+j}x_2^{1-\a+i-k}x_3^{-j}x_4^{(4\b+2-\a+2i-2j+2k)}\rt],\no\\
&&Ch_{KCoker_{\F_{(\a;\b)}}}(q,x_1,x_2,x_3,x_4)
=\frac{q^{\frac{1}{2}\a(2\b+1)}}{\prod_{n=1}^{\infty}(1-q^n)^6}
\sum_{l_1,l_2=1}^{\infty}(-1)^{l_1+l_2}q^{\frac{1}{2}(l_1^2+l_2^2-
(1-2\a)l_1+l_2)}\no\\
&&~~~~~~~~~~~~~~\times \sum_{i,j,k\in
\Z}\lt[q^{\frac{1}{2}(i^2+j^2+k^2-2kj+(1-2l_1)i+(1+2l_1)j+(1+2l_2)k)}\rt.\no\\
&&~~~~~~~~~~~~~~\times
 \lt.x_1^{\a+j}x_2^{1-\a+i-k}x_3^{-j}x_4^{(4\b+2-\a+2i-2j+2k)}\rt],\no\\
&&Ch_{KKer_{\F_{(\a;\b)}}}(q,x_1,x_2,x_3,x_4)
=\frac{q^{\frac{1}{2}\a(2\b+1)}}{\prod_{n=1}^{\infty}(1-q^n)^6}
\sum_{l_1,l_2=1}^{\infty}(-1)^{l_1+l_2}q^{\frac{1}{2}(l_1^2+l_2^2+
(1-2\a)l_1+l_2)}\no\\
&&~~~~~~~~~~~~~~\times \sum_{i,j,k\in
\Z}\lt[q^{\frac{1}{2}(i^2+j^2+k^2-2kj+(1+2l_1)i+(1-2l_1)j+(1+2l_2)k)}\rt.\no\\
&&~~~~~~~~~~~~~~~\times
 \lt.x_1^{\a+j}x_2^{1-\a+i-k}x_3^{-j}x_4^{(4\b+2-\a+2i-2j+2k)}\rt],\no
\eea
\noindent and 
\bea
&&Sch_{KKer_{\F_{(\a;\b)}}}(q,x_1,x_2,x_3,x_4)
=\frac{q^{\frac{1}{2}\a(2\b+1)}}{\prod_{n=1}^{\infty}(1-q^n)^6}
\sum_{l_1,l_2=1}^{\infty}(-1)^{l_1+l_2}q^{\frac{1}{2}(l_1^2+l_2^2-
(1-2\a)l_1-l_2)}\no\\
&&~~~~~~~~~~~~~~~\times \sum_{i,j,k\in
\Z}\lt[ 
(-1)^{1-\a+i-j+k}q^{\frac{1}{2}(i^2+j^2+k^2-2kj+(1-2l_1)i+(1+2l_1)
j+(1-2l_2)k)}\rt.\no\\
&&~~~~~~~~~~~~~~~~\times\lt.
x_1^{\a+j}x_2^{1-\a+i-k}x_3^{-j}x_4^{(4\b+2-\a+2i-2j+2k)}
\rt],\no\\
&&Sch_{CKer_{\F_{(\a;\b)}}}(q,x_1,x_2,x_3,x_4)
=\frac{q^{\frac{1}{2}\a(2\b+1)}}{\prod_{n=1}^{\infty}(1-q^n)^6}
\sum_{l_1,l_2=1}^{\infty}(-1)^{l_1+l_2}q^{\frac{1}{2}(l_1^2+l_2^2+
(1-2\a)l_1-l_2)}\no\\
&&~~~~~~~~~~~~~~~~\times \sum_{i,j,k\in
\Z}\lt[(-1)^{1-\a+i-j+k}q^{\frac{1}{2}(i^2+j^2+k^2-2kj+(1+2l_1)i+(1-2l_1)j
+(1-2l_2)k)}\rt.\no\\
&&~~~~~~~~~~~~~~~~\times\lt.
x_1^{\a+j}x_2^{1-\a+i-k}x_3^{-j}x_4^{(4\b+2-\a+2i-2j+2k)}\rt],\no\\
&&Sch_{KCoker_{\F_{(\a;\b)}}}(q,x_1,x_2,x_3,x_4)
=\frac{q^{\frac{1}{2}\a(2\b+1)}}{\prod_{n=1}^{\infty}(1-q^n)^6}
\sum_{l_1,l_2=1}^{\infty}(-1)^{l_1+l_2}q^{\frac{1}{2}(l_1^2+l_2^2-
(1-2\a)l_1+l_2)}\no\\
&&~~~~~~~~~~~~~~~~\times \sum_{i,j,k\in
\Z}\lt[(-1)^{1-\a+i-j+k}q^{\frac{1}{2}(i^2+j^2+k^2-2kj+(1-2l_1)i
+(1+2l_1)j+(1+2l_2)k)}\rt.\no\\
&&~~~~~~~~~~~~~~~~\times\lt.
x_1^{\a+j}x_2^{1-\a+i-k}x_3^{-j}x_4^{(4\b+2-\a+2i-2j+2k)}\rt],\no\\
&&Sch_{CCoker_{\F_{(\a;\b)}}}(q,x_1,x_2,x_3,x_4)
=\frac{q^{\frac{1}{2}\a(2\b+1)}}{\prod_{n=1}^{\infty}(1-q^n)^6}
\sum_{l_1,l_2=1}^{\infty}(-1)^{l_1+l_2}q^{\frac{1}{2}(l_1^2+l_2^2+
(1-2\a)l_1+l_2)}\no\\
&&~~~~~~~~~~~~~~~~\times \sum_{i,j,k\in
\Z}\left[(-1)^{1-\a+i-j+k}q^{\frac{1}{2}(i^2+j^2+k^2-2kj+(1+2l_1)i
+(1-2l_1)j+(1+2l_2)k)}\rt.\no\\
&&~~~~~~~~~~~~~~~~\times\lt.
x_1^{\a+j}x_2^{1-\a+i-k}x_3^{-j}x_4^{(4\b+2-\a+2i-2j+2k)}\rt].\no
\eea
\end{Proposition}
\noindent{\it Proof}. We sketch the proof of these three propositions.
Since  $q^{-d}x_1^{h_1}x_2^{h_2}x_3^{h_3}x_4^{h_4}$ and 
 $(-1)^{N_f}q^{-d}x_1^{h_1}x_2^{h_2}x_3^{h_3}x_4^{h_4}$ commute with 
the BRST charges $Q^{(i)}_{l}$ and $[Q^{(1)}_l,Q^{(2)}_{l'}]=0$,
 the trace over $Ker$ and $Coker$ can be 
written as the sum of trace over $\F_{(*;\b)}^{(l_1,l_2)}$. 
The latter can be 
computed by the technique introduced in \cite{Cla73}.

\vspace{0.5truecm}

Note that $\F_{(*;\b)}^{(1,1)}=\F_{(*;\b-1)}$, we have 
\begin{Corollary}
: The following relations hold for any $\a\in \Z$ and $\b$,
\bea
&&Ch_{CCoker_{\F_{(*;\b+1)}}}=Ch_{KKer_{\F_{(*,\b)}}},\\
&&Sch_{CCoker_{\F_{(*;\b+1)}}}=Sch_{KKer_{\F_{(*,\b)}}}.
\eea
\end{Corollary}

Now, we study the $\uqg2h$-module structures of $\F_{(*;\b)}$ for 
$*=(0,1),(1,0)$, $\a\in \Z$. Set 
\begin{eqnarray*}
&&\l^0_{\b}=\L_0+4\b\L_4,~~~~~ \l^3_{\b}=\L_3+(4\b+3)\L_4,\\
&&\l_{\a,\b}=\a\L_1+(1-\a)\L_2+(4\b+2-\a)\L_4, ~~ {\rm for }\; \a\in\Z, 
\end{eqnarray*}
\noindent for arbitrary $\b$  and 
\begin{eqnarray*}
&&|\l^0_{\b}>=|\b,-\b,\b,-\b;0,0>\in \F_{((0,1);\b)},\\
&&|\l^3_{\b}>=|\b+1,-\b-1,\b+1,-\b;0,0>\in \F_{((1,0);\b)},\\
&&|\l_{\a,\b}>=|\b+1,-\b-1+\a,\b,-\b;-\a,0>\in \F_{(\a;\b)},\; \a\in\Z.
\end{eqnarray*}
\noindent The above vectors play the role of the highest weight vectors of 
 $\uqg2h$-modules in the Fack spaces $\F_{(*;\b)}$ with $*=(0,1),(1,0)$,
$\a\in\Z$. One can verify that ($l=1,2$ below)
\bea 
&&\eta^l_0|\l^0_{\b}>=\eta^l_0|\l^3_{\b}>=\eta^l_0|\l_{\a,\b}>=0,~~
 {\rm for} ~\a=0,-1,-2,\cdots,\\
&&\eta^2_0|\l_{\a,\b}>=0,~~~~  \eta^1_0|\l_{\a,\b}>\neq 0,~~ 
 {\rm for} ~\a=1,2,\cdots.
\eea
\noindent It follows that the modules 
$KKer_{\F_{((0,1);\b)}}$, $KKer_{\F_{((1,0);\b)}}$ 
are highest weight $\uqg2h$-modules with highest weights $\l^0_{\b}$, 
$\l^3_{\b}$ respectively, while 
 $KKer_{\F_{(\a;\b)}}$ $(\a=0,-1,-2,\cdots)$ and 
 $CKer_{\F_{(\a;\b)}}$ $(\a=1,2,\cdots)$ 
 are highest weight $\uqg2h$-modules with highest 
weights $\l_{\a,\b}$.
We denote them by $\V(\l^0_{\b})$, $\V(\l^3_{\b})$ and $\V(\l_{\a,\b})$, 
respectively:  
\begin{eqnarray*}
\V(\l^0_{\b})&\cong & KKer_{\F_{((0,1);\b)}},~~~~~
\V(\l^3_{\b})\cong  KKer_{\F_{((1,0);\b)}},\\
\V(\l_{\a,\b}) &\cong & KKer_{\F_{(\a;\b)}},~~ {\rm
for}~\a=0,-1,-2,\cdots\\
&\cong & CKer_{\F_{(\a;\b)}},~~{\rm for}~
\a=1,2,3,\cdots.
\end{eqnarray*}
\noindent It is expected that the modules $\V(\l^0_{\b})$
, $\V(\l^0_{\b})$, $\V(\l_{\a,\b})$ are also irreducible with 
respect to the action of $\uqg2h$. Namely,

\vspace{0.4truecm}
\noindent{\bf Conjecture 2}: {\it
$\V(\l^0_{\b})$, $\V(\l^3_{\b})$ and $\V(\l_{\a,\b})$ are the 
irreducible highest weight $\uqg2h$-modules with the highest weight 
 $\l^0_{\b}$, $\l^3_{\b}$ and $\l_{\a,\b}$, respectively, i.e. 
\bea
&&\V(\l_{\a,\b})=V(\l_{\a,\b}),~~ \a\in\Z,\\
&&\V(\l^0_{\b})=V(\l^0_{\b}),~~~~~ \V(\l^3_{\b})=V(\l^3_{\b}).
\eea}

\subsection{Vertex operators of $\uqg2h$ }
In this subsection we study the action of the type I and type II vertex
operators  of $\uqg2h$ 
  on the highest weight $\uqg2h$-modules.

Using the bosonic representations of the vertex operators
 (\ref{Vertex-operator}), we have the homomorphisms of 
$\uqg2h$-modules:
\bea
\phi(z)&:&~\F_{(\a;\b)}\longrightarrow\F_{(\a-1;\b)}\otimes V_z,\no\\
\psi(z)&:&~\F_{(\a;\b)}\longrightarrow V_z\otimes\F_{(\a-1;\b)},
   \label{Inertwining-relation1}\\
\phi^*(z)&:&~\F_{(\a;\b)}\longrightarrow\F_{(\a+1;\b)}\otimes V^{*S}_z,\no\\
\psi^*(z)&:&~\F_{(\a;\b)}\longrightarrow V^{*S}_z\otimes\F_{(\a+1;\b)}.
   \label{Intertwining-relation2}
\eea
Then we  consider the vertex operators which intertwine the highest 
weight $\uqg2h$-modules by using the above results. For $\a\not\in \Z$,
we have by  conjecture 1: 

\vspace{0.4truecm}
\noindent{\bf Conjecture 3}: 
{\it The following vertex operators 
associated with the level-one irreducible highest weight 
$\uqg2h$-modules exist:
\bea 
\Phi(z)^{\l_{\a-1,\b}~V}_{\l_{\a,\b}}(z)&:&~~
V(\l_{\a,\b})\longrightarrow V(\l_{\a-1,\b})\otimes
V_z,\no\\
\Psi(z)_{\l_{\a,\b}}^{V~\l_{\a-1,\b}}&:&~~V(\l_{\a,\b})\longrightarrow
V_z\otimes V(\l_{\a-1,\b}),\no\\
\Phi(z)^{\l_{\a+1,\b}~V^*}_{\l_{\a,\b}}(z)&:&~~
V(\l_{\a,\b})\longrightarrow V(\l_{\a+1,\b})\otimes
V^{*S}_z,\no\\
\Psi(z)_{\l_{\a,\b}}^{V^*~\l_{\a+1,\b}}&:&~~V(\l_{\a,\b})\longrightarrow
V^{*S}_z\otimes V(\l_{\a+1,\b}),
\eea
\noindent for $\a \not\in \Z$.}

\vspace{0.4truecm}
It is easy to see that the bosonized vertex operators 
(\ref{Vertex-operator}) also commute (or anti-commute) with $\eta^1_0$ 
and $\eta^2_0$. Noting this property, the homomorphisms
(\ref{Inertwining-relation1})-(\ref{Intertwining-relation2}) and  
conjecture 2, we have 

\vspace{0.4truecm}
\noindent{\bf Conjecture 4}: {\it  For $\a\in\Z$, 
the following vertex operators
associated with the level-one irreducible highest weight
$\uqg2h$-modules exist:
\bea
\Phi(z)^{\l_{\a-1,\b}~V}_{\l_{\a,\b}}(z)&:&~~
V(\l_{\a,\b})\longrightarrow V(\l_{\a-1,\b})\otimes
V_z,~~~ \a\neq 1,\no\\
\Psi(z)_{\l_{\a,\b}}^{V~\l_{\a-1,\b}}&:&~~V(\l_{\a,\b})\longrightarrow
V_z\otimes V(\l_{\a-1,\b}),~~~ \a\neq 1,\no\\
\Phi(z)^{\l_{\a+1,\b}~V^*}_{\l_{\a,\b}}(z)&:&~~
V(\l_{\a,\b})\longrightarrow V(\l_{\a+1,\b})\otimes
V^{*S}_z,~~~ \a\neq 0,\no\\
\Psi(z)_{\l_{\a,\b}}^{V^*~\l_{\a+1,\b}}&:&~~V(\l_{\a,\b})\longrightarrow
V^{*S}_z\otimes V(\l_{\a+1,\b}),~~~ \a\neq 0.
\eea}
\vspace{0.4truecm}
\vskip.3in
\subsection*{ Acknowledgements.}

This work has been  financially supported by Australian Research 
Council large, small and QEII fellowship grants. We would like to thank
Prof. B.Y. Hou for his encouragement and useful discussions.
 W.L. Yang  thanks Y.-Z. Zhang and  the department 
of Mathematics, the University of Queensland,
for their kind  hospitality. W.L. Yang was also partially supported by the
National Natural Science Foundation of China.

\subsection*{Appendix A}
In this appendix, we give the normal ordered relations of the fundamental 
bosonic fields:

\begin{eqnarray*}
&&:e^{H^i(z;\b_1)}::e^{H^j(w;\b_2)}:=(z-wq^{\b_1+\b_2})^{a_{ij}}
:e^{H^i(z;\b_1)+H^j(w;\b_2)}:~,\\
&&:e^{H^i(z;\b_1)}::e^{H^{*j}(w;\b_2)}:=(z-wq^{\b_1+\b_2})^{\delta_{ij}}
:e^{H^i(z;\b_1)+H^{*j}(w;\b_2)}:~,\\      
&&:e^{H^{*i}(z;\b_1)}::e^{H^{j}(w;\b_2)}:=(z-wq^{\b_1+\b_2})^{\delta_{ij}}
:e^{H^{*i}(z;\b_1)+H^{j}(w;\b_2)}:~,\\
&&:e^{H^{*1}(z;\b_1)}::e^{H^{*2N-1}(w;\b_2)}:=
:e^{H^{*1}(z;\b_1)+H^{*2N-1}(w;\b_2)}:~,\\
&&:e^{B_{2N}(z;\b_1)}::e^{H^{i}(w;\b_2)}:=:e^{H^{i}(w;\b_2)}:
:e^{B_{2N}(z;\b_1)}:=:e^{B_{2N}(z;\b_1)+H^{i}(w;\b_2)}:~,\\
&&:e^{B_{1}(z;\b_1)}::e^{H^{i}(w;\b_2)}:=:e^{H^{i}(w;\b_2)}: 
:e^{B_{1}(z;\b_1)}:=:e^{B_{1}(z;\b_1)+H^{i}(w;\b_2)}:~,\\ 
&&:e^{B_{2N}(z;\b_1)}::e^{B_{2N}(w;\b_2)}:
=:e^{B_{2N}(w;\b_2)}::e^{B_{2N}(z;\b_1)}:
=:e^{B_{2N}(z;\b_1)+B_{2N}(w;\b_2)}:~,\\
&&:e^{B_{1}(z;\b_1)}::e^{B_{1}(w;\b_2)}:
=:e^{B_{1}(w;\b_2)}::e^{B_{1}(z;\b_1)}:
=:e^{B_{1}(z;\b_1)+B_{1}(w;\b_2)}:~,\\
&&:e^{B_{2N}(z;\b_1)}::e^{H^{*2l-1}(w;\b_2)}:
=z^{\frac{1}{2N}}(1-\frac{w}{z}q^{\b_1+\b_2})
:e^{B_{2N}(z;\b_1)+H^{*2l-1}(w;\b_2)}:~,\\
&&:e^{B_{2N}(z;\b_1)}::e^{H^{*2l}(w;\b_2)}:
=:e^{B_{2N}(z;\b_1)+H^{*2l}(w;\b_2)}:~,\\
&&:e^{c^l(z)}::e^{c^{l'}(z)}:=(z-w)^{\delta_{ll'}}:e^{c^l(z)+c^{l'}(z)}:~,
\end{eqnarray*}
\noindent where $i,j=1,2,\cdots,2N-1$ and  $l,l'=1,2,\cdots,N$.

\subsection*{Appendix B}
By using   Theorem 1,
the integral
expressions of the bosonized vertex operators (\ref{Vertex-operator}) and 
the technique  in \cite{AJMP}, one can check the following
relations

\begin{itemize}

\item For the type I vertex operators:

\begin{eqnarray*}
&&[\phi_k(z),f_l]=0 ~~{\rm if}~k\neq l,l+1,~~~~
[\phi_l(z),f_l]_{q^{(-1)^l}}=0,\\
&&[\phi_{l+1}(z),f_l]_{q^{(-1)^l}}=
(-1)^l\phi_l(z),\\
&&[\phi_k(z),e_l]=0 ~~{\rm if}~k\neq l,~~~~
[\phi_{l}(z),e_l]=-q^{h_l}\phi_{l+1}(z),\\
&&q^{h_l}\phi_l(z)q^{-h_l}=q^{(-1)^l}\phi_l(z),\\   
&&q^{h_l}\phi_k(z)q^{-h_l}=\phi_k(z)~{\rm if}~k\neq l,l+1,~~~~
q^{h_l}\phi_{l+1}(z)q^{-h_l}=q^{(-1)^l}\phi_{l+1}(z),\\
\\
&&[\phi^*_k(z),f_l]=0 ~~{\rm if}~k\neq l,l+1,~~~~
[\phi^*_{l+1}(z),f_l]_{q^{(-1)^{l+1}}}=0,\\
&&[\phi^*_k(z),e_l]=0 ~~{\rm if}~k\neq l+1,~~~~
[\phi^*_{l+1}(z),e_l]=(-1)^lq^{h_l+(-1)^l}\phi^*_{l}(z),\\
&&[\phi^*_{l}(z),f_l]_{q^{(-1)^{l+1}}}=q^{(-1)^{l+1}}\phi^*_{l+1}(z),~~~~
q^{h_l}\phi^*_l(z)q^{-h_l}=q^{(-1)^{l+1}}\phi^*_l(z),\\
&&q^{h_l}\phi^*_k(z)q^{-h_l}=\phi^*_k(z)~{\rm if}~k\neq l,l+1,~~~~
q^{h_l}\phi^*_{l+1}(z)q^{-h_l}=q^{(-1)^{l+1}}\phi^*_{l+1}(z).
\\
\end{eqnarray*}
\item For the type II vertex operators:
\begin{eqnarray*}
&&[\psi_k(z),e_l]=0 ~~{\rm if}~k\neq l,l+1,~~~~
[\psi_{l+1}(z),e_l]_{q^{(-1)^{l+1}}}=0,\\ 
&&[\psi_{l}(z),e_l]_{q^{(-1)^{l+1}}}=\psi_{l+1}(z),\\
&&[\psi_k(z),f_l]=0 ~{\rm if}~k\neq l+1,~~~~ 
[\psi_{l+1}(z),f_l]=(-1)^{l+1}q^{-h_l}\psi_{l}(z),\\
&&q^{h_l}\psi_l(z)q^{-h_l}=q^{(-1)^l}\psi_l(z),~~~~ 
q^{h_l}\psi_{l+1}(z)q^{-h_l}=q^{(-1)^l}\psi_{l+1}(z),\\
&&q^{h_l}\psi_k(z)q^{-h_l}=\psi_k(z)~~{\rm if}~k\neq l,l+1,\\
\\
&&[\psi^*_k(z),e_l]=0 ~~{\rm if}~k\neq l,l+1,~~~~ 
[\psi^*_{l}(z),e_l]_{q^{(-1)^l}}=0,\\
&&[\psi^*_k(z),f_l]=0 ~~{\rm if}~k\neq l,~~~~    
[\psi^*_{l}(z),f_l]=-q^{-h_l-(-1)^l}\psi^*_{l+1}(z),\\
&&[\psi^*_{l+1}(z),e_l]_{q^{(-1)^{l}}}=(-1)^{l+1}q^{(-1)^l}
\psi^*_{l}(z),~~~~ 
q^{h_l}\psi^*_l(z)q^{-h_l}=q^{(-1)^{l+1}}\psi^*_l(z),\\
&&q^{h_l}\psi^*_k(z)q^{-h_l}=\psi^*_k(z)~{\rm if}~k\neq l,l+1,~~~~
q^{h_l}\psi^*_{l+1}(z)q^{-h_l}=q^{(-1)^{l+1}}\psi^*_{l+1}(z).
\end{eqnarray*}
\end{itemize}


\end{document}